\documentclass{article}
\usepackage{amssymb, amsmath}
\usepackage[latin9]{inputenc}
\usepackage[francais]{babel}
\usepackage[affil-it]{authblk}
\usepackage{graphicx}
\usepackage{hyperref}
\usepackage{color}
\usepackage[onehalfspacing]{setspace}
\usepackage{ulem}
\makeatletter
\def\imod#1{\allowbreak\mkern10mu({\operator@font mod}\,\,#1)}
\makeatother
\usepackage{fullpage}
\usepackage{fancyhdr}

\begin{document}

\title{La fonte algébrique \\ des \textit{Méthodes nouvelles de la mécanique céleste} \\ d'Henri Poincaré}
\author{Fr{\'e}d{\'e}ric Brechenmacher
\thanks{Electronic address: \texttt{frederic.brechenmacher@euler.univ-artois.fr \\ \textit{Ce travail a  b\'{e}n\'{e}fici\'{e} d'une aide de l'Agence Nationale de la Recherche : projet CaaF\'{E} (ANR-10-JCJC 0101)}}}} 
\affil{Universit{\'e} d'Artois \\ Laboratoire de math{\'e}matiques de Lens (EA 2462) \\
rue Jean Souvraz S.P. 18, 62307 Lens Cedex France.
 \\ \& \\
 {\'E}cole polytechnique \\ D{\'e}partement humanités et sciences sociales \\
91128 Palaiseau Cedex, France. \\}
 \date{}
\maketitle

\begin{abstract}
L'approche de Poincaré sur le problème des trois corps a souvent été célébrée comme un point d'origine de l'étude des systèmes dynamiques et de la théorie du chaos. Cet article propose de porter un autre regard sur cette approche en analysant le rôle qu'y jouent   des pratiques algébriques spécifiques de manipulation des systèmes linéaires.  Ces pratiques algébriques permettent de mieux  cerner la spécificité de l'approche de Poincaré en restituant les temporalités et dimensions collectives d'une oeuvre souvent célébrée pour son individualité et son caractère de point d'origine. 

Bien qu'elles soient passées inaperçues de l'historiographie, ces pratiques jouent un véritable rôle de modèle pour  \textit{Les méthodes nouvelles de la mécanique céleste} qu'elles soutiennent comme par une sorte de fonte algébrique. Comme un alliage conserve la trace des métaux qui le constituent, cette fonte témoigne de cadres collectifs dans lesquels s'inscrit Poincaré. Mais tout comme la structure de fonte d'un immeuble disparaît derrière l'ornementation créative d'une façade, le moule algébrique de la stratégie de Poincaré se brise en engendrant  \textit{Les méthodes nouvelles de la mécanique céleste}. 

Une version éditée de cette prépublication sera publiée prochainement dans le journal \textit{L'astronomie} sous le titre "L'approche de Poincaré sur le problème des trois corps". Les  principaux arguments donnés dans cet article, ainsi que les enjeux historiographiques et  les aspects mathématiques brièvement évoqués  aux paragraphes 4 et 5,  seront par ailleurs développés dans une publication prochaine - "The algebraic cast of Poincaré's \textit{Méthodes nouvelles}" - dont cette prépublication constitue un résumé.
\bigskip

\textbf{Abstracted version of  "The algebraic cast of Poincaré's \textit{Méthodes nouvelles}"}.
Poincaré's approach to the three body problem has often been celebrated as a starting point of chaos theory in relation to the investigation of dynamical systems. Yet, Poincaré's strategy can also be analyzed as molded on - or casted in - some specific algebraic practices for manipulating systems of linear equations. These practices shed new light on both the novelty and the collective dimensions of Poincaré's \textit{Les méthodes nouvelles de la mécanique céleste}. As the structure of a cast-iron building  may be less noticeable than its creative façade, the algebraic cast of Poincaré's strategy is broken out of the mold in generating the novel methods of celestial mechanics. But as the various components that are mixed in some casting process can still be detected in the resulting alloy, the algebraic cast of the \textit{Méthodes nouvelles}  points to some collective dimensions of Poincaré's methods.

An edited version of the present preprint is to be published in the journal \textit{L'astronomie} under the title "L'approche de Poincaré sur le problème des trois corps".  This publication is  an abstract in French language of a forthcoming paper -  "The algebraic cast of Poincaré's \textit{Méthodes nouvelles}"  - which will  develop its main claims as well as the historiographical and mathematical issues raised in section 4 and section 5.

\end{abstract}

\pagestyle{fancy}
\fancyhead{}
\lhead{La fonte algébrique des \textit{Méthodes nouvelles de la mécanique céleste}}
\rhead{ F. BRECHENMACHER}
 \renewcommand{\headsep}{18pt}

\tableofcontents

\section{Introduction}

\begin{quote}
Les personnes qui s'intéressent aux progrès de la Mécanique Céleste, mais qui ne peuvent les suivre que de loin, doivent éprouver quelque étonnement en voyant combien de fois on a démontré la stabilité du système solaire.
Lagrange l'a établie d'abord, Poisson l'a démontrée de nouveau, d'autres démonstrations sont venues depuis, d'autres viendront encore. [...] L'étonnement de ces personnes redoublerait sans doute, si on leur disait qu'un jour peut-être un mathématicien fera voir, par un raisonnement rigoureux, que le système planétaire est instable. [...] On peut démontrer, dans certains cas particuliers, que les éléments de l'orbite d'une planète redeviendront une infinité de fois très voisins des éléments initiaux, et cela est probablement vrai aussi dans le cas général, mais cela ne suffit pas ; il faudrait faire voir que non seulement ces éléments finiront par reprendre leurs valeurs primitives, mais qu'ils ne s'en écarteront jamais beaucoup.

Cette dernière démonstration, on ne l'a jamais donnée d'une manière rigoureuse et il est probable que la proposition n'est pas rigoureusement vraie. Ce qui est vrai seulement, c'est que les éléments [de l'orbite] ne pourront s'écarter sensiblement de leur valeur primitive qu'avec une extrême lenteur et au bout d'un temps tout à fait énorme. Aller plus loin, affirmer que ces éléments resteront, non pas très longtemps, mais toujours, compris entre des limites étroites, c'est ce que nous ne pouvons faire.\cite{Poincar1898}
\end{quote}

À partir des années 1950, \textit{Les méthodes nouvelles de la mécanique céleste} ont souvent été présentées comme fondatrices de la théorie du chaos. Différents aspects des travaux d'Henri Poincaré ont été mis en avant selon les acteurs impliqués dans cette théorie en plein essor : l'étude qualitative des équations différentielles, envisagée dès 1881 en lien avec les trajectoires géométriques des corps célestes ; l'attention à la stabilité globale d'ensembles de telles trajectoires ; la notion de bifurcation, liée à l'étude de systèmes différentiels variables en fonction d'un paramètre ; ou encore l'introduction de concepts probabilistes en dynamique. Le célèbre théorème de récurrence de Poincaré énonce ainsi qu'un système mécanique isolé retourne à un état proche de son état initial sauf pour un ensemble de cas de probabilité nulle :
\begin{quote}
Il y a donc une infinité de solutions particulières qui sont instables, au sens que nous venons de donner à ce mot et une infinité d'autres qui sont stables. J'ajouterai que les premières sont exceptionnelles (ce qui permet de dire qu'il a stabilité en général). Voici ce que j'entends par là, car ce mot par lui même n'a aucun sens. Je veux dire qu'il y a une probabilité nulle pour que les conditions initiales du mouvement soient celles qui correspondent à une solution instable.\cite{Poincar1891}
\end{quote}
De fait, les trois tomes des \textit{Méthodes nouvelles} ont suscité des lectures diverses qui ne mettent pas toutes l'accent sur les mêmes résultats, approches ou concepts.\cite{AubinDahan} Mais ces lectures rétrospectives ont aussi relégué au second plan certains aspects auparavant considérés comme importants. L'un d'entre eux se manifeste dans l'usage que fait Poincaré des trajectoires périodiques :\footnote{Notons que c'est à la classification des solutions périodiques que le mathématicien Jacques Hadamard attribuait l'une des principales nouveautés de l'approche de Poincaré peu après le décès de ce dernier.\cite[p.643]{Hadamard1913}}
\begin{quote}
Ce sont celles où les distances des trois corps sont des fonctions périodiques du temps ; à des intervalles périodiques, les trois corps se retrouvent donc dans les mêmes positions relatives.\cite{Poincar1891}
\end{quote}
La classification de ces solutions particulières supporte l'étude des trajectoires plus complexes, comme les solutions asymptotiques qui, sans être périodiques, tendent à le devenir au bout d'un temps infini ; et comme les solutions  doublement asymptotiques - ou homoclines - qui s'enroulent autour de trajectoires périodiques.

Les trajectoires périodiques ne valent donc pas pour elles-mêmes mais pour l'usage qu'en fait Poincaré : elles permettent l'introduction de systèmes \textit{linéaires} et, par là, de pratiques \textit{algébriques} spécifiques. Bien qu'elles soient passées inaperçues de l'historiographie, ces pratiques jouent un véritable rôle de modèle pour les \textit{Méthodes nouvelles}. Elles sont non seulement importantes pour saisir la nouveauté de ces méthodes mais permettent également de restituer les temporalités et dimensions collectives d'une oeuvre souvent célébrée pour son individualité et son caractère de point d'origine. De fait, l'usage de systèmes linéaires en mécanique porte la marque des grands traités de la fin du XVIII\up{e} siècle. Il semble donc à première vue paradoxal que Poincaré leur associe la nouveauté de ses \textit{Méthodes} :
\begin{quote}
L'étude des inégalités séculaires par le moyen d'un système d'équations différentielles linéaires à coefficients constants peut donc être regardée comme se rattachant plutôt aux méthodes nouvelles qu'aux méthodes anciennes.\cite{Poincar1892}
\end{quote}
Nous verrons que Poincaré s'appuie non seulement sur d'anciens travaux de Joseph-Louis Lagrange et Pierre-Simon Laplace mais aussi sur des développements qui ont été donnés à ces travaux tout au long du XIX\up{e} siècle, en un va-et-vient entre astronomie, géométrie, arithmétique, algèbre et analyse. Nous aurons ainsi l'occasion de porter un regard original sur les relations entre la mécanique céleste et d'autres branches des sciences mathématiques. 

\section{Le problème des trois corps}
Selon la loi de Newton, l'attraction mutuelle entre les planètes perturbe l'ellipse qu'adopterait une planète unique en orbite autour du soleil. Bien que faibles, ces inégalités séculaires s'opposent au calcul des éphémérides à longue échéance. Plus encore, elles posent la question de la validité et de l'universalité de loi de la gravitation :
\begin{quote}
Il ne s'agit pas seulement, en effet, de calculer les éphémérides des astres quelques années d'avance pour les besoins de la navigation ou pour que les astronomes puissent retrouver les petites planètes déjà connues. Le but final de la Mécanique céleste est plus élevé ; il s'agit de résoudre cette importante question : la loi de Newton peut-elle expliquer à elle seule tous les phénomènes astronomiques ? [...] elle a pour expression mathématique une équation différentielle, et pour obtenir les coordonnées des astres, il faut intégrer cette équation. [...] Quel sera le mouvement de $n$ points matériels, s'attirant mutuellement en raison directe de leurs masses et en raison inverse du carré des distances ? Si $n=2$, c'est à dire si l'on a affaire à une planète isolée et au Soleil, en négligeant les perturbations dues aux autres planètes, l'intégration est facile ; les deux corps décrivent des ellipses, en se conformant aux lois de Kepler. La difficulté commence si le nombre $n$ des corps est égal à trois : le problème des trois corps a défié jusqu'ici tous les efforts des analystes.\cite{Poincar1891}
\end{quote}
Sur le long terme, les inégalités séculaires  peuvent en effet  engendrer des perturbations importantes des orbites et menacer ainsi la stabilité du système solaire :
 \begin{quote}
Une des questions qui ont le plus préoccupé les chercheurs est celle de la stabilité du système solaire. C'est à vrai dire une question mathématique plutôt que physique. Si l'on découvrait une démonstration générale et rigoureuse, on n'en devrait pas conclure que le système solaire est éternel. Il peut en effet être soumis à d'autres forces que celle de Newton, et les astres ne se réduisent pas à des points matériels. [...] on n'est pas absolument certain qu'il n'existe pas de milieu résistant ; d'autre part les marées absorbent de l'énergie qui est incessamment convertie en chaleur par la viscosité des mers , et cette énergie ne peut être empruntée qu'à la force vive des corps célestes. [...] Mais toutes ces causes de destruction agiraient beaucoup plus lentement que les perturbations, et si ces dernières n'étaient pas capables d'en altérer la stabilité, le système solaire serait assuré d'une existence beaucoup plus longue.\cite{Poincar1891}
 \end{quote}

Rappelons que le problème des trois corps a été proposé pour le grand prix organisé à l'occasion du soixantième anniversaire du roi Oscar II de Suède et Norvège.\cite{Barrow-Greene1996}. Bien que lauréat, le mémoire soumis par Poincaré présente des conclusions erronées quant à la stabilité des trajectoires. C'est en rectifiant cette erreur que Poincaré a introduit les solutions homoclines,\cite{Anderson1994} classiquement considérées comme la première description d'un comportement chaotique. Un mémoire modifié a été publié en 1890,\cite{Poincar1890} avant que Poincaré ne s'attelle à la rédaction des \textit{Méthodes nouvelles} de 1892 à 1899. 

Un corps céleste y est assimilé à une masse ponctuelle de coordonnées $(x_1,..., x_n)$ décrivant, en fonction du temps  $t$,  une trajectoire  exprimée par des fonctions analytiques $X_i$  des coordonnées selon le système différentiel :
\[
\frac{dx_i}{dt}=X_i \ (i=1,..., n)
\]
Ces équations ne peuvent être résolues de manière exacte. Il s'agit d'en approcher les solutions générales par des solutions particulières :
\begin{quote}
Le mouvement des trois astres dépend en effet de leurs positions et de leurs vitesses initiales. Si l'on se donne ces conditions initiales du mouvement, on aura défini une solution particulière du mouvement. [...] Il peut arriver que les orbites des trois corps se réduisent à des ellipses. La position et la vitesse initiales de notre satellite auraient pu être telles que la Lune fût constamment pleine ; elles auraient pu être telles que la Lune fût constamment nouvelle [...] les phases auraient pu suivre des lois bien étranges ; dans une des solutions possibles, la Lune, d'abord nouvelle, commence par croître ; mais, avant d'atteindre le premier quartier, elle se met à décroître pour redevenir nouvelle et ainsi de suite ; elle aura donc constamment la forme d'un croissant.\cite{Poincar1891}
\end{quote}
Si certaines de ces solutions particulières ne sont "intéressantes que par leur bizarrerie", d'autres sont susceptibles d'applications astronomiques. Dans les années 1880, l'astronome George William Hill s'est notamment appuyé sur des solutions périodiques dans ses travaux sur la théorie de la lune.\cite{Hill1877} \cite{Hill1886} Poincaré s'en inspire, bien que son objectif ne soit pas d'étudier ces solutions pour elles-mêmes :
\begin{quote}
En effet, il y a une probabilité nulle pour que les conditions initiales du mouvement soient précisément celles qui correspondent à une solution périodique. Mais il peut arriver qu'elles en diffèrent très peu, et cela a lieu justement dans les cas où les méthodes anciennes ne sont plus applicables. On peut alors avec avantage prendre la solution périodique comme première approximation. [...] Il y a même plus : voici un fait que je n'ai pu démontrer rigoureusement mais qui me paraît pourtant très vraisemblable [...] on peut toujours trouver une solution périodique (dont la période peut, il est vrai, être très longue), telle que la différence entre les deux solutions soit aussi petite qu'on le veut, pendant un temps aussi long qu'on le veut. D'ailleurs, ce qui nous rend ces solutions périodiques si précieuses, c'est qu'elles sont, pour ainsi dire, la seule brèche par où nous puissions essayer de pénétrer dans une place jusqu'ici réputée inabordable.\cite{Poincar1891}
\end{quote}
Les trajectoires périodiques donnent lieu à deux types de méthodes d'approximations. La \textit{méthode de variation} consiste à étudier des ensembles de solutions d'un \textit{même} système différentiel. Une solution périodique donnée a-t-elle un comportement similaire à celui d'une autre solution dont les conditions initiales sont proches ? La \textit{méthode de perturbation} consiste, quant à elle, à faire varier le système différentiel en fonction d'un petit paramètre $\mu$. Dans le "problème restreint" des trois corps qu'aborde Poincaré, cette approche est légitimée par la petitesse des masses des planètes comparées à celle du soleil 
\begin{quote}
Ce cas est celui du problème où l'on suppose que les trois corps se meuvent dans un même plan, que la masse du troisième est nulle, que les deux premiers décrivent des circonférences concentriques autour de leur centre de gravité commun.\cite{Poincar1891}
\end{quote}
Le troisième corps est hors d'état de troubler les deux autres corps; il s'agit d'en comprendre la trajectoire en fonction du rapport $\mu$ des masses des deux autres corps. Lorsque $\mu=0$, cette trajectoire est une ellipse képlérienne. Que se passe-t-il lorsque $\mu$ n'est pas nul, mais très petit ? 
\begin{quote}
Avons-nous le droit [de] conclure [qu'un système admettant des solutions périodiques pour  $\mu = 0$] en admettra encore pour les petites valeurs de $\mu$ ? [...] La première solution périodique qui ait été signalée pour le cas où $\mu>0$ est celle qu'a découverte Lagrange et où les trois corps décrivent des ellipses képleriennes semblables, pendant que leurs distances mutuelles restant dans un rapport constant. [...] M. Hill, dans ses très remarquables recherches sur la théorie de la Lune en a étudié une autre [...]. Je montrerai [...] comment on peut prendre une solution périodique comme point de départ d'une série d'approximations successives, et étudier ainsi les solutions qui en diffèrent fort peu.\cite[p. 106 \& 153]{Poincar1892}
\end{quote}
Une grande partie des travaux de Poincaré est consacrée à  démontrer  l'existence de solutions périodiques pour certaines conditions initiales et à l'analyse de leurs comportements après perturbation. Comme nous allons le voir, cette stratégie d'approximations par des trajectoires périodiques prend modèle sur des pratiques algébriques développées au XVIII\up{e} siècle pour mathématiser les petites oscillations de systèmes de corps.

\section{Des petites oscillations d'une corde à celles des trajectoires périodiques}

A la fin du XVIII\up{e} siècle, les inégalités séculaires ont été étudiées sur le modèle de la mathématisation donnée à certains problèmes de cordes vibrantes. Mentionnons Jean le Rond d'Alembert qui, en 1743, a étudié les petites oscillations $\xi_i(t)$ d'une corde lestée de deux masses en négligeant les termes non linéaires dans les développements en séries des équations de la dynamique. En 1766, Lagrange a généralisé cette approche aux oscillations de $n$ masses, mathématisées par un système de n équations différentielles linéaires à coefficients constants :\cite{Lagrange1766}
\[
\frac{d\xi_i}{dt}=\sum_{j=1,n} A_{i,j}\xi_j
\]
L'intégration du système est obtenue par sa décomposition en $n$ équations indépendantes. Cette méthode mathématise la propriété mécanique selon laquelle les oscillations d'une corde lestée de $n$ masses peuvent se décomposer en oscillations indépendantes de $n$ cordes lestées d'une seule masse. Si l'on pose comme inconnue $S$, la période d'une telle oscillation propre, on est amené au calcul du déterminant\footnote{Cette équation correspond à l'équation caractéristique de la matrice symétrique $A=(A_{ij})$ qu'elle vise à diagonaliser en déterminant ses valeurs propres. Cette perspective ne s'est cependant imposée que dans les années 1930. Elle introduit implicitement des conceptions anachroniques qui masquent d'anciennes formes d'organisations du savoir qui ne se traduisent pas dans un cadre disciplinaire comme celui de l'algèbre linéaire.}
\[
\begin{vmatrix}
A_{1,1}-S & A_{1,2} & ... & A_{1,n} \\ 
A_{2,1}  & A_{2,2}-S & ... & A_{2,n} \\ 
... &...& ... & ... \\ 
A_{n,1} & A_{n,2} & ... & A_{n,n} -S \\ 
\end{vmatrix}
\]
qui donne une équation algébrique de degré $n$, souvent désignée sous le nom d'"équation en $S$". À chaque racine $\alpha_i$ de cette équation est associée une oscillation propre $\xi_i(t)=e^{\alpha_i t}$. Dans le cas où les $n$ racines sont distinctes, on obtient ainsi $n$ solutions indépendantes à partir desquelles on peut exprimer toutes les solutions du système. 

Dans les années 1780, Lagrange et Laplace ont transféré cette mathématisation au cas des oscillations des planètes sur leurs orbites.\cite{Lagrange1774}  \cite{Laplace1775} Pour cette raison, l'équation en $S$ a par la suite également été désignée sous le nom d'"équation à l'aide de laquelle on détermine les inégalités séculaires des planètes" (que nous abrégerons par "équation séculaire").

Chez Poincaré, la méthode d'intégration de Lagrange est sous-jacente à la stratégie d'approximations par des solutions périodiques.

\begin{quote}
Il y a peu de chances pour que, dans aucune application, les conditions initiales du mouvement soient exactement celles qui correspondent à une solution périodique ; mais il peut arriver qu'elles en diffèrent fort peu. Si alors on considère les coordonnés des trois corps dans leur mouvement véritable, et, d'autre part, les coordonnées qu'auraient ces trois mêmes corps dans la solution périodique, la différence reste très petite au moins pendant un certain temps et l'on peut, dans une première approximation, négliger le carré de cette différence.\cite[p.162]{Poincar1892}
\end{quote}

Une petite variation $\xi_i(t)$ d'une trajectoire périodique $\phi_i(t)$ est ainsi mathématisée par le système linéaire (2). Si l l'équation en $S$ a toutes ses racines distinctes, alors 
\[
\xi_i(t)=C_1e^{\alpha_1t}\lambda_{1,i}(t)+C_2e^{\alpha_2t}\lambda_{2,i}(t)+...+C_ne^{\alpha_nt}\lambda_{n,1}(t)
\]
où les $\lambda_{i,j} (t)$ sont des séries trigonométriques uniformément convergentes de même période que $\phi_i(t)$.

Les coefficients $\alpha_i$ sont dénommés "exposants caractéristiques" par Poincaré. Observons que ces derniers mêlent des significations mécaniques et algébriques qui ne sont pas identiques chez Lagrange et Poincaré. Nous avons vu que la méthode du premier s'appuie sur une représentation mécanique des racines $\alpha_i$ comme les périodes des petites oscillations propres de l'orbite elliptique des planètes. Or, chez Poincaré, ces petites oscillations opèrent cette fois sur une trajectoire périodique, dans le voisinage de laquelle elles engendrent un faisceau d'autres trajectoires. Comme nous allons le voir, le comportement de ces faisceaux est contrôlé par la nature algébrique - et notamment la multiplicité - des exposants caractéristiques.

\section{Stabilité et multiplicité des exposants caractéristiques}

\begin{quote}
Le mot stabilité a été entendu sous les sens les plus différents, et la différence de ces divers sens deviendra manifeste si l'on se rappelle l'histoire de la Science.

Lagrange a démontré qu'en négligeant les carrés des masses, les grands axes des orbites demeurent invariables. Il voulait dire par là qu'avec ce degré d'approximation les grands axes peuvent se développer en séries dont les termes sont de la forme $Asin(\alpha t+ \beta)$, $A$, $\alpha$ et $\beta$ étant des constantes.

Il en résulte que, si ces séries sont uniformément convergentes, les grands axes demeurent compris entre certaines limites ; le système des astres ne peut donc pas passer par toutes les situations compatibles avec les intégrales des forces vives et des aires, et plus il repassera une infinité de fois aussi près que l'on voudra de sa situation initiale. C'est la stabilité complète. 

Poussant plus loin l'approximation, Poisson a annoncé ensuite que la stabilité subsiste quand on tient compte des carrés des masses et qu'on en néglige les cubes. Mais cela n'avait pas le même sens. Il voulait dire que les grands axes peuvent se développer en séries contenant non seulement des termes de la forme $Asin(\alpha t+\beta)$ mais des termes de la forme $Atsin(\alpha t+\beta)$. La valeur du grand axe éprouve alors de continuelles oscillations, mais rien ne prouve que l'amplitude de ces oscillations ne crois pas indéfiniment avec le temps. Nous pouvons affirmer que le système repassera toujours une infinité de fois aussi près qu'on voudra de sa situation initiale mais non qu'il ne s'en éloignera pas beaucoup. Le mot  de stabilité n'as donc pas le même sens pour Lagrange et pour Poisson.

[...] Pour qu'il y ait stabilité complète dans le problème des trois corps, il faut trois conditions :
\begin{enumerate}
\item Qu'aucun des trois corps ne puisse s'éloigner indéfiniment
\item Que deux des corps ne puissent se choquer et que la distance de ces deux corps ne puisse descendre au-dessous d'une certaine limite
\item Que le système vienne repasser une infinité de fois aussi près que l'on veut de sa situation initiale.
\end{enumerate}
[...] Il y a un cas [celui du problème restreint] où, depuis longtemps, on a démontré que la première condition est remplie. Nous allons voir que la troisième l'est également. Quant à la deuxième, je ne puis rien dire.\cite{Poincar1892} 
\end{quote}
Au XVIII\up{e} siècle, un critère de stabilité mécanique a été donné en fonction des racines de l'équation en $S$. d'Alembert a ainsi énoncé qu'une corde infiniment peu dérangée d'un état d'équilibre tend à y revenir si les $\alpha_i$ sont \textit{réels, négatifs et distincts}. En effet, les solutions prennent dans ce cas une forme périodique, $sin(\alpha_i t)$, tandis que des racines imaginaires impliquent des oscillations exponentielles. L'occurrence de racines multiples, $\alpha_i=\alpha_j$, pose un problème plus délicat : la méthode d'intégration semble alors mise en échec puisque celle-ci repose sur une décomposition en $n$ équations indépendantes, chacune associée à une racine distincte. Dans ce cas, les solutions ont et exprimées sous la forme $tsin(\alpha_it)$ : "le temps sort du sinus" et génère des oscillations non périodiques et non bornées. 

Dans le cas de la mécanique céleste traité par Lagrange, de telles expressions ont été désignées sous le nom de \textit{termes séculaires}. Pour les petites oscillations d'une corde, la stabilité était une donnée de la situation mécanique étudiée, tandis que les inégalités séculaires  des planètes posent le problème de la stabilité du système solaire. Le transfert d'une même mathématisation des oscillations d'une corde à celles des planètes a ainsi engendré des discussions sur la nature algébrique des racines de l'équation séculaire. La célèbre démonstration par Laplace de la stabilité du système du monde (dans le cas de l'approximation linéaire) contient ainsi une preuve que ces racines sont réelles. Cette preuve s'appuie sur la symétrie des systèmes différentiels de la mécanique - propriété mise en évidence par les travaux de Lagrange sur l'équation séculaire. Elle conclut à la disparition de tout terme séculaire.

À partir des travaux de Siméon Denis Poisson en 1808,\cite{Poisson1809} la question de la stabilité a été traitée par la prise en compte de termes non linéaires dans les développements en séries des coordonnées des astres. En 1856, Urbain Le Verrier a notamment montré que ces termes ne permettent pas seulement des approximations plus fines des trajectoires mais peuvent induire des altérations importantes des orbites.\footnote{Pour une histoire sur le temps long du problème de la stabilité du système solaire de Newton à Poincaré, voir \cite{Laskar1992}}\cite{leverrier1856} Au cours du XIX\up{e} siècle, de nombreux astronomes, , parmi lesquels Newcomb, Linstedt, Gyldén, Delaunay ou Hill, ont proposé de nouvelles séries dans l'objectif d'en éliminer tout terme séculaire. La notion de stabilité a alors progressivement été associée avec l'obtention de développements en séries strictement trigonométriques dont les premiers termes décroissent rapidement. Il en a découlé une conception de la "convergence astronomique" des séries distincte de la notion mathématique de convergence vers une limite finie. 

Or, Poincaré a démontré que les séries des astronomes sont divergentes au sens mathématique :
\begin{quote}
[...] les méthodes de M. Gyldén et celles de M. Lindstedt ne donnent en effet, si loin que l'on pousse l'approximation, que des termes périodiques, de sorte que tous les éléments des orbites ne peuvent éprouver que des oscillation autour de leur valeur moyenne. La question serait donc résolue, si ces développement étaient convergents. Nous savons malheureusement qu'il n'en est rien.\cite{Poincar1891}
\end{quote} 
Poincaré a ainsi été amené à poser la question de la stabilité différemment de la majorité des travaux du XIX\up{e} siècle.\footnote{Sur les différentes notions de stabilité chez Poincaré et leurs circulations, voir \cite{Roque2011}}. Il définit notamment la "stabilité au sens de Lagrange" d'une solution périodique en revenant au cas linéaire et à la discussion sur les racines de l'équation en $S$, c'est à dire les exposants caractéristiques qu'il appelle aussi "coefficients de stabilité". Les raisonnements du XVIII\up{e} sont repris presque mot pour mot, ainsi que leur conclusion :
\begin{quote}
En résumé, $\xi_i$ peut dans tous les cas être représenté par une série toujours convergente. Dans cette série, le temps peut entrer sous le signe sinus ou cosinus, ou par l'exponentielle $e^{\alpha t}$, ou enfin en dehors des signes trigonométriques ou exponentiels. Si tous les coefficients de stabilité sont réels, négatifs et distincts, le temps n'apparaîtra que sous les signes sinus et cosinus et il y aura stabilité temporaire. Si l'un des coefficients est positif ou imaginaire, le temps apparaîtra sous un signe exponentiel ; si deux des coefficients sont égaux ou que l'un deux soit nul, le temps apparaît en dehors de tout signe trigonométrique ou exponentiel.  Si donc tous les coefficients ne sont pas réels, négatifs et distincts, il n'y a pas en général de stabilité temporaire. [...] Il ne faut pas toutefois entendre ce mot de stabilité au sens absolu. En effet, nous avons négligé les carrés des $\xi$ [...] et rien ne prouve qu'en tenant compte de ces carrés, le résultat ne serait pas changé. Mais nous pouvons dire au moins que les $\xi$ [...], s'ils sont originairement très petits, resteront très petits pendant très longtemps; Nous pouvons exprimer ce fait en disant que la solution périodique jouit, sinon de la stabilité séculaire, du moins de la stabilité temporaire.\cite{Poincar1892}
\end{quote}
Ce critère de stabilité est un premier pas vers l'analyse de flots de trajectoires dans le voisinage d'une solution périodique : si cette solution reste stable, les trajectoires approchées lui restent proches tandis que des solutions périodiques instables permettent d'introduire des trajectoires plus complexes comme les solutions asymptotiques. 

\section{La culture algébrique portée par l'équation séculaire}

Nous avons jusqu'à présent éclairé les rôles joués par les systèmes linéaires chez Poincaré en mettant ce dernier en relation directe avec Lagrange. Mais cette relation ne saurait être considéré comme exclusive. Si les références des Méthodes nouvelles aux grands traités du XVIII\up{e} siècle sont fréquentes, la lecture en est faite au prisme de nombreux autres travaux publiés au XIX\up{e} siècle. Il nous faut à présent changer d'échelle d'analyse en portant notre attention à la longue durée et aux dimensions collectives de l'approche de Poincaré. 

Comme nous l'avons déjà évoqué, le grand traité de Poincaré a marqué une date dans l'histoire de l'astronomie. Comme les traités de Lagrange ou Laplace, il a, à son tour, fait l'objet de lectures multiples et fragmentées. Au XX\up{e} siècle, Poincaré a notamment souvent été présenté comme le premier à rouvrir le dossier de la stabilité du système du monde après Laplace. Mais si la stabilité a pu être considéré comme acquise au XIX\up{e} siècle, les méthodes développées en mécanique céleste n'en ont pas moins suscité de nombreux travaux dans d'autres domaines. 

On trouve ainsi au XIX\up{e} siècle de très nombreuses références à l'équation séculaire. Celles-ci manifestent rarement des préoccupations en astronomie mais identifient une culture algébrique partagée à l'échelle européenne.\cite{Brechenmacher:2007b}. Les procédés élaborés par Lagrange pour manipuler des systèmes linéaires par la décomposition polynomiale de l'équation séculaire sont au coeur de cette culture commune ; nous avons vu qu'ils irriguent encore en profondeur les Méthodes nouvelles. Mais cette culture algébrique ne se limite pas à la technicité de certains procédés opératoires. Ces derniers véhiculent également un idéal de généralité qui se manifeste notamment dans l'ambition de Poincaré de traiter de situations à $n$ variables. Cette généralité participe du caractère spécial de l'équation séculaire : l'équation étant de degré $n$, elle ne peut en général être résolue par radicaux mais la nature réelle de ses racines peut être déduite de la propriété de symétrie des systèmes mécaniques qui l'engendrent. 

Tout au long du XIX\up{e} siècle, cette équation spéciale a supporté des analogies et des transferts de procédés opératoires entre différentes branches des sciences mathématiques. Nous avons déjà évoqué la représentation mécanique sous-jacente aux procédés de décompositions de Lagrange. Plus tard, en 1829, Augustin Louis Cauchy a interprété ces procédés comme revenant à  la recherche des axes principaux d'une conique. Plus tard encore, ces mêmes procédés ont servi de modèle à des théories de l'élasticité et de la lumière, ils ont été intégrés à l'analyse complexe, à la théorie des formes quadratiques etc.

Mais l'équation séculaire a surtout véhiculé des problématiques spécifiques relatives aux racines multiples. Ces dernières donnent lieu à des difficultés aussi délicates que celles qui ont amené à l'introduction des nombres complexes. C'est d'ailleurs en partie pour surmonter le problème de la multiplicité des racines de l'équation séculaire que Cauchy a élaboré son calcul des résidus en analyse complexe. À partir des années 1850, plusieurs approches algébriques du problème de la multiplicité ont été développées de manières distinctes : introduction par James Joseph Sylvester de la notion de matrice en géométrie analytique, théorie des formes quadratiques de Charles Hermite, diviseurs élémentaires de Karl Weierstrass, réduction canonique de Camille Jordan en théorie des groupes. 

Le rôle de modèle que joue une certaine pratique des systèmes linéaires dans les Méthodes nouvelles s'inscrit donc dans une culture algébrique largement partagée au XIX\up{e} siècle. La manière dont Poincaré aborde la question de l'existence de solutions périodiques en témoigne. Cette question est en effet traitée par une discussion sur la multiplicité des racines de l'équation en $S$ comme l'illustre notamment l'énoncé suivant : "les solutions périodiques du problème des trois corps ont deux exposants caractéristiques nuls, mais elles n'en ont que deux". 

Mais Poincaré a aussi développé une approche qui lui est propre en s'appuyant sur deux lignes distinctes de développements des travaux sur l'équation séculaire. Les méthodes de Jordan en théorie des groupes lui permettent de traiter les occurrences de racines multiples dans des systèmes linéaires qui, contrairement à ceux de Lagrange, ne sont pas symétriques et ne peuvent en général être décomposés en $n$ équations indépendantes. Poincaré s'est approprié ces méthodes en les mêlant à l'héritage des travaux d'Hermite.\cite{Brechenmacher:2012a}. Cet héritage se manifeste dès 1881 lorsque Poincaré légitime son approche qualitative des équations différentielles par analogie avec le théorème de Sturm qui donne le nombre de racines réelles d'une équation algébrique :
\begin{quote}
Dans les cas élémentaires, l'expression des inconnues par les symboles usuels fournit en général aisément à leur égard tous les renseignement que l'on se propose d'obtenir. [...] Pour peu que la question se complique, il en est autrement et, tout au moins, on peut dire que la lecture, si j'ose m'exprimer ainsi, faite par le mathématicien des documents qu'il possède, comporte deux grandes étapes, l'une que l'on peut appeler qualitative, l'autre quantitative. Ainsi, par exemple, pour étudier une équation algébrique, on commence par rechercher, à l'aide du théorème de Sturm, quel est le nombre des racines réelles : c'est la partie qualitative ; puis on calcule la valeur numérique de ces racines, ce qui constitue l'étude quantitative de l'équation. [...] C'est naturellement par la partie qualitative qu'on doit aborder la théorie de toute fonction et c'est pourquoi le problème qui se présente en premier lieu est le suivant : \textit{Construire les courbes définies par des équations différentielles}.\cite{Poincar1881f}
\end{quote}
En effet, dans les années 1850, Hermite a cherché à obtenir une nouvelle démonstration du théorème de Sturm en explorant le cas particulier de l'équation séculaire.\cite[p.124-132]{Sinaceur1991} Sa méthode permet d'analyser la variation du nombre de racines distinctes de cette équation en fonction de l'inconnue $S$.\footnote{Cette approche est basée sur la loi d'intertie des formes quadratiques, ; pour davantage de développements sur les parties 4. et 5. de cet article, voir  \cite{Brechenmacher:2013}} Elle donne un modèle pour aborder le problème crucial de la perturbation des solutions périodiques en fonction du paramètre $\mu$. Poincaré transfère ainsi la notion de multiplicité des racines des équations algébriques aux trajectoires des systèmes différentiels. Il énonce par exemple que "les solutions périodiques disparaissent par couples à la façon des racines réelles des équations algébriques". Il s'agit là d'un résultat clé pour analyser l'effet de perturbations sur des solutions périodiques stables : nous avons vu que l'occurrence d'une racine multiple dans l'équation en $S$ est une cause d'instabilité ; de même, si une trajectoire "périodique perd la stabilité ou l'acquiert", c'est "qu'elle se sera confondue avec une autre solution périodique, avec laquelle elle aura échangé sa stabilité".

\section{Conclusion}

Nous avons vu dans cet article quelques exemples de la manière dont l'architecture des \textit{Méthodes nouvelles }est soutenue par une sorte de \textit{fonte algébrique}. Comme un alliage conserve la trace des métaux qui le constituent, cette fonte témoigne de cadres collectifs dans lesquels s'inscrit l'approche de Poincaré : culture algébrique portée par l'équation séculaire sur le temps long ; héritages, sur des temps plus courts, des travaux d'Hermite et de Jordan. Mais tout comme la structure de fonte d'un immeuble disparaît derrière l'ornementation créative d'une façade, le moule algébrique de la stratégie de Poincaré se brise en engendrant les nouvelles méthodes de mécanique céleste.

Cette situation éclaire d'autres aspects dont nous ne pouvons traiter ici, comme les procédés itératifs de la célèbre méthode des sections qui témoignent eux-aussi de l'héritage de travaux d'Hermite sur le théorème de Sturm, notamment développé dans les années 1870-1880 par Edmond Laguerre.\footnote{Voir à ce sujet \cite{Brechenmacher:2013}} Les rapports entre la mécanique céleste et les autres branches des sciences mathématiques au XIX\up{e} siècle sont ainsi loin de se réduire à un va-et-vient entre application et abstraction. Non seulement des procédés algébriques spécifiques ont émergé de travaux de mécanique céleste, mais l'équation séculaire a joué un rôle de référence partagée tout au long du XIX\up{e} siècle. Cette équation a supporté le transfert de procédés opératoires entre différents domaines, circulation qui a enrichi ces procédés de significations nouvelles avant que Poincaré ne les investisse à nouveau en mécanique céleste. Ce constat ne manque pas d'analogie avec celui qu'énonce Poincaré à propos de la complexité des solutions des équations différentielles : des trajectoires peuvent fortement s'éloigner des conditions initiales mais néanmoins revenir dans leur voisinage après une longue durée.

\bibliography{bibliographie}
\bibliographystyle{apalike-fr}

\end{document}